\documentclass[11pt]{amsart}

\usepackage{amsmath}
\usepackage{amssymb}
\usepackage{bbm}
\usepackage{pdfsync}

\newcommand{\C}{{\mathbb C}}       
\newcommand{\R}{{\mathbb R}}       

\newcommand{\HH}{{\mathcal H}}

\newcommand{\UU}{{\mathcal U}}

\newcommand{\dist}{{\rm dist}}

\newcommand{\rf}[1]{{(\ref{#1})}}

\newcommand{\supp}{\operatorname{supp}}

\newcommand{\ve}{{\varepsilon}}
\newcommand{\vv}{{\vspace{2mm}}}
\newcommand{\vvv}{{\vspace{3mm}}}
\newcommand{\wt}[1]{{\widetilde{#1}}}

\newcommand{\pv}{\operatorname{pv}}

\newtheorem{theorem}{Theorem}[section]

\newtheorem*{theorem*}{Theorem}

\theoremstyle{definition}

\theoremstyle{remark}
\newtheorem{rem}[theorem]{Remark}

\numberwithin{equation}{section}

\newcommand{\brem}{\begin{rem}}
\newcommand{\erem}{\end{rem}}

\begin{document}

\title[Jump formulas on rectifiable sets]{Jump formulas for singular integrals and layer potentials on rectifiable sets}

\author{Xavier Tolsa}

\address{Xavier Tolsa
\\
ICREA, Passeig Llu\'{\i}s Companys 23 08010 Barcelona, Catalonia\\
 Departament de Matem\`atiques, and BGSMath
\\
Universitat Aut\`onoma de Barcelona
\\
08193 Bellaterra (Barcelona), Catalonia.
}
\email{xtolsa@mat.uab.cat} 

\thanks{Partially supported by 2017-SGR-0395 (Catalonia) and MTM-2016-77635-P (MINECO, Spain).
}

\begin{abstract}
In this paper the jump formulas for the double layer potential and other
singular integrals are proved for arbitrary rectifiable sets, by defining suitable non-tangential limits.
The arguments are quite straightforward and only require some Calder\'on-Zygmund techniques.
\end{abstract}

\maketitle

\section{Introduction}

Let $K:\R^{n+1}\to \R$ (or $K:\R^{n+1}\to \R^{n+1}$) be an odd  Calder\'on-Zygmund kernel of homogeneity $-n$. That is, $K$ satisfies the estimate
\begin{equation}\label{eq1}
|K(x)|\leq \frac{ C}{|x|^n}\quad \mbox{ for $x\neq 0$,}
\end{equation}
and there exists some $\eta>0$ such that
\begin{equation}\label{eq2}
|K(x) - K(x+y)|\leq C\,\frac{|y|^\eta}{|x|^{n+\eta}}\quad \mbox{ if $|y|\leq \frac12|x|$,}
\end{equation}
and moreover $K$ is of the form $K(x) = \dfrac{\Omega(x)}{|x|^n}$,
with $\Omega: \R^{n+1}\to \R^{n+1}$ being homogeneous of degree $0$.
Given a signed Radon measure $\nu$ in $\R^{n+1}$ we denote
$$T\nu(x) = \int K(x-y)\,d\nu(y),$$
whenever the integral makes sense,
and for $\ve>0$
$$T_\ve\nu(x) = \int_{|x-y|>\ve} K(x-y)\,d\nu(y)$$
and also
$$T_*\nu(x) =\sup_{\ve>0} |T_\ve\nu(x)\quad \mbox{ and }\quad\pv T\nu(x) = \lim_{\ve\to0} T_\ve \nu(x),$$
whenever the last limit exists.

Given a positive Radon measure $\mu$, we consider the operators $T_\mu$, $T_{\mu,\ve}$, $T_{\mu,*}$ defined by
$$T_\mu f = T(f\mu),\qquad T_{\mu,\ve} f = T_\ve(f\mu),\qquad  T_{\mu,*} f = T_*(f\mu),$$
for $f\in L^1_{loc}(\mu)$. As usual, one says that $T_\mu$ is bounded in $L^p(\mu)$ if the operators $T_{\mu,\ve}$ are bounded in $L^p(\mu)$ uniformly on $\ve>0$. Analogously, $T$ is bounded from the space of (signed) finite
Radon measures $M(\R^{n+1})$ to $L^{1,\infty}(\mu)$ if the operators $T_\ve$ are bounded 
from $M(\R^{n+1})$ to $L^{1,\infty}(\mu)$ uniformly on $\ve>0$. That is, there exists some constant $C$ such that 
for all $\nu\in M(\R^{n+1})$, all $\ve>0$, and all $\lambda>0$,
$$\mu\big(\big\{x\in\R^{n+1}:|T_\ve\nu(x)|>\lambda\big\}\big)\leq C\,\frac{\|\nu\|}\lambda.$$
Recall that if $\mu$ satisfies the polynomial growth condition
$$\mu(B(x,r))\leq C\,r^n\quad \mbox{ for all $x\in\R^{n+1}$, $r>0$,}$$
then the $L^2(\mu)$ boundedness of $T_\mu$ implies the $L^p(\mu)$ boundedness of $T_\mu$ and $T_{\mu,*}$
for $1<p<\infty$, and also the boundedness from $M(\R^{n+1})$ to $L^{1,\infty}(\mu)$, by results essentially due to Nazarov, Treil and Volberg in the non-doubling setting (see Chapter 2 of
\cite{Tolsa-llibre}, for example).

A set $E\subset\R^{n+1}$ is called $n$-rectifiable if there is a collection of $n$-dimensional Lipschitz
 (possibly rotated) graphs $\Gamma_k\subset\R^{n+1}$ such that 
 $$\HH^n\Big(E\setminus \textstyle\bigcup_{k\geq1}\Gamma_k\Big)=0.$$
 An equivalent definition (that will be used below) is obtained if the graphs are required to be $C^1$ and
 with maximal slope as close to $0$ as wished.

As shown by Mas \cite{Mas}, for any $n$-rectifiable set $E\subset\R^{n+1}$ and any measure $\nu\in M(\R^{n+1})$, 
if the  kernel $K$ is odd and $C^2$ away from the origin and satisfies
\begin{equation}\label{eq4}
|\nabla^j K(x)|\leq \frac C{|x|^{n+j}}\quad \mbox{ for $j=0,1,2$ and $x\neq 0$},
\end{equation}
then $\pv T\nu(x)$ exists for $\HH^n$-a.e.\ $x\in E$. Previously, under the same assumptions in \cite{Tolsa-PLMS} it was proved that if $\mu$ is uniformly $n$-rectifiable, then $T_\mu$ is bounded in $L^2(\mu)$.
Under stronger smoothness conditions on the kernel, this had been proved previously  by David and Semmes
\cite{DS1}.

In this paper we will see how one can prove a very general form of the jump formulas (also known as Plemelj formulas in the case of the Cauchy transform) for the operators described above. Recall that these formulas play an important role in the
solution of the Dirichlet problem by the method of layer potentials.
To state precisely our results we need to introduce some additional notation. For a point $x\in\R^{n+1}$, a unit vector $u$, 
and an aperture parameter $a\in(0,1)$ we consider the one sided cone with axis in the direction of $u$ defined by
$$X_a(x,u)=\bigl\{y\in\R^{n+1}:(y-x)\cdot u> a|y-x|\bigr\}.$$
We say that $E$ has a classical tangent $n$-plane at $x\in E$ if there exists a unit vector $u$ such that, for all
$a\in(0,1)$, there exists some $r_0>0$ such that
$$E\cap \bigl(X_a(x,u) \cup X_a(x,-u)\bigr) \cap B(x,r) =\varnothing \quad\mbox{ for $0<r\leq r_0$.}$$
The $n$-plane $L$ orthogonal to $u$ through $x$ is called a (classical) tangent $n$-plane at $x$.

We say that $E$ has an approximate tangent $n$-plane at $x\in E$ if there exists a unit vector $u$ such that, for all
$a\in(0,1)$,
$$\lim_{r\to0} \frac{\HH^n\bigl(E\cap \bigl(X_a(x,u) \cup X_a(x,-u)\bigr)\cap B(x,r)\bigr)}{r^n}=0.$$
The $n$-plane $L$ orthogonal to $u$ through $x$ is called approximate tangent $n$-plane.

If $\HH^n(E)<\infty$ (or $\HH^n|_E$ is locally finite) and $E$ is $n$-rectifiable, it is well known that there is a unique
approximate tangent $n$-plane at $\HH^n$-a.e.\ $x\in E$. We denote by $L_x$ the approximate tangent $n$-plane at $x$
and by $N_x$  a unit vector orthogonal to $L_x$. We also write
$$X_a^+(x) = X_a(x,N_x),\qquad X_a^-(x) = X_a(x,-N_x),\qquad X_a(x) = X_a^+(x) \cup X_a^-(x).$$
As mentioned above, the $n$-plane $L_x$ is uniquely defined for $\HH^n$-a.e.\ $x\in E$. On the other hand, 
for $\HH^n$-a.e.\ $x\in E$ there are two possible choices for $N_x$, depending on the sense of the normal to $L_x$ that one chooses. For us, the choice $x\mapsto N_x$ does not matter as soon as it is $\HH^n$-measurable (or Borel, say if $E$ is Borel). 
Notice that any $y\in X_a^+(x)\cup X_a^-(x)$ satisfies
$$\dist(y,L_x)> a\,|y-x|>0.$$
Fix $b\in(0,a)$, so that $\bar B(y,b|y-x|)\cap L_x=\varnothing$ for all $y\in X_a^+(x)\cup X_a^-(x)$.
We define the non-tangential limits
\begin{equation}\label{eqdtt}
T^+\nu(x) = \lim_{X_a^+(x)\ni y\to x} T_{b|x-y|}\nu(y),\qquad T^-\nu(x) = \lim_{X_a^-(x)\ni y\to x} T_{b|x-y|}\nu(y),
\end{equation}
whenever they exist. Note that we use the truncated operators $T_{b|x-y|}$ in these definitions, which may appear rather unusual.

Our main result is the following.

\begin{theorem}\label{teo1}
Let $T$ be the operator associated with an odd Calder\'on-Zygmund kernel of homogeneity $-n$ satisfying
\rf{eq4}, $C^2$ away from the origin. Let $E\subset \R^{n+1}$ be an $n$-rectifiable set and let $\nu\in M(\R^{n+1})$. For fixed $a\in
(0,1)$ and $b\in(0,a)$ as above, the non-tangential limits 
$T^+\nu(x)$ and $T^-\nu(x)$ exist $\HH^n$-a.e.\ $x\in E$ and moreover the following identities hold
$\HH^n$-a.e.\ $x\in E$ too:
\begin{equation}\label{eqga}
 \frac12 \bigl(T^+\nu(x) + T^-\nu(x)\bigr)= \pv T\nu(x),
 \end{equation}
 and
\begin{equation}\label{eqgb}
\frac12 \bigl(T^+\nu(x) - T^-\nu(x)\bigr) = C_K(N_x)\,f(x),
\end{equation}
where $f$ is the density 
$$f(x) = \frac{d\nu}{d\HH^n|_E}(x),$$
and $C_K(N_x)$ is defined by
\begin{equation}\label{eqckn}
C_K(N_x)= \int_{L(N_x)} \frac{\Omega(y+N_x)- \Omega(y-N_x)}{2(|y|^2 + 1)^{n/2}}\,d\HH^n(y),
\end{equation}
where $L(N_x)$ is the hyperplane orthogonal to $N_x$ through the origin.
\end{theorem}
\vv

Some remarks are in order:
\begin{itemize}
\item Write $\nu = \nu_{abs} + \nu_s$, where $\nu_{abs}$ is absolutely continuous with respect to 
$\HH^n|_E$ and $\nu_s$ is singular with respect to $\HH^n|_E$. Then, for $\HH^n$-a.e.\ $x\in E$,
$$\frac{d\nu}{d\HH^n|_E}(x) = \frac{d\nu_{abs}}{d\HH^n|_E}(x).$$
\vv

\item If $E$ has a classical tangent $n$-plane at $x$ and $\supp\nu\subset E$, then
\begin{equation}\label{eqtt5}
T^+\nu(x) = \lim_{X_a^+(x)\ni y\to x} T\nu(y),\qquad T^-\nu(x) = \lim_{X_a^-(x)\ni y\to x} T\nu(y),
\end{equation}
because, for any $r>0$ small enough, $\bar B(y,b|x-y|)\cap E = \varnothing$ for all $y\in X_a(x)\cap B(x,r)$ and thus $T_{b|x-y|}\nu(y)= T\nu(y)$.
The identities in \rf{eqtt5} may fail (even at $\HH^n$-a.e.\ $x\in E$) if there is not a classical tangent.
It is not difficult to construct an example showing this.

Notice that if $E$ is $n$-AD regular and rectifiable, then any approximate tangent $n$-plane is also a classical tangent $n$-plane.\vv

\item Theorem \ref{teo1} is also valid for any operator $T$ associated with a kernel satisfying 
\rf{eq1} and \rf{eq2}, assuming that for any $n$-dimensional Lipschitz graph $\Gamma$ (possibly rotated)
the operator $T_{\HH^n|_\Gamma}$ is bounded in $L^2(\HH^n|_\Gamma)$ and that $\pv T\HH^n|_{\Gamma\cap E}$ exists $\HH^n$-a.e.\ in $\Gamma$ for all $E\subset \Gamma$ with $\HH^n(E)<\infty$. 

\vv

\item If $T_{\HH^n|_E}$ is bounded in $L^2(\HH^n|_E)$, then the jump formulas \rf{eqga} and \rf{eqgb} can be extended
to any $f\in L^p(\HH^n|_E)$, $1\leq p<\infty$, by very standard arguments, just by replacing $\nu$ by $f\,\HH^n|_E$.\vv

\item Theorem \ref{teo1} can be extended to $n$-rectifiable sets $E\subset \R^d$ with $n<d-1$. In this case, given an approximate tangent $n$-plane $L_x$ at $x\in E$, there are infinitely many unit vectors orthogonal to $L_x$. Then we just fix one such vector $N_x$, we define $X_a(x,N_x)$ as above and we denote
$$X_a(x,N_x,L_x) = X_a(x,N_x) \cap [L_x,N_x],$$
where $[L_x,N_x]$ is the affine $(n+1)$-plane that contains $L_x$ and $x+N_x$.
We set
$$X_a^+(x) = X_a(x,N_x,L_x),\qquad X_a^-(x) = X_a(x,-N_x,L_x),$$
and define $T^+\nu(x)$, $T^-\nu(x)$ as in \rf{eqdtt}.
Then Theorem \ref{teo1} still holds, with the same proof essentially, with the integral in
\rf{eqckn} over the $n$-plane parallel to $L_x$ through the origin.
\vv

\item When $K$ is the  $n$-dimensional Riesz kernel, i.e.,
$$K(x) = \frac x{|x|^{n+1}},\qquad \Omega(x) = \frac x{|x|},$$
for all $y\in L(N_x)$ we have 
$$\Omega(y+N_x)- \Omega(y-N_x) = \frac{y+N_x}{|y+N_x|} - \frac{y-N_x}{|y-N_x|}
= \frac{2N_x}{(|y|^2+1)^{1/2}}$$
and thus
\begin{align*}
C_K(N_x) & = \int_{L(N_x)} \frac{N_x}{(|y|^2 + 1)^{(n+1)/2}}\,d\HH^n(y) \\
&= N_x\int_{\R^n} \frac{1}{(|y|^2 + 1)^{(n+1)/2}}\,d\HH^n(y) = \frac{\pi^{(n+1)/2}}{\Gamma\bigl(\tfrac{n+1}2\bigr)}\,N_x = \frac{\omega_{n}}2\,N_x,
\end{align*}
where we assumed $\HH^n$ to be defined so that it coincides with $n$-dimensional Lebesgue measure in $\R^n$
and $\omega_{n}$ is the $n$-dimensional volume of the unit sphere in $\R^{n+1}$.

So in the case of the double layer potential, defined by
$$R_{\HH^n|_E} f(x) = \frac1{\omega_n} \int \frac{x-y}{|x-y|^{n+1}} \,N_y\,f(y)\,d\HH^n|_E,$$
by applying \rf{eqga} and \rf{eqgb} to the vectorial measure $f(y)\,N_y\,d\HH^n|_E(y)$, it follows that
$$
 \frac12 \bigl(R_{\HH^n|_E}^+f(x) + R_{\HH^n|_E}^-f(x)\bigr)= \pv R_{\HH^n|_E}f(x),
$$
 and
$$
\frac12 \bigl(R_{\HH^n|_E}^+f(x) - R_{\HH^n|_E}^-f(x)\bigr) = \frac{1}2\,N_x\cdot N_x\,f(x) = \frac12\,f(x)
$$
for $\HH^n$-a.e.\ $x\in E$.
\vv

\item For an odd integer $j\geq1$, consider the kernel $K$ in the complex plane defined by
$$K(z) = \frac{z^j}{|z|^{j+1}},\qquad \Omega(z) = \frac{z^j}{|z|^j},$$
for $z\in\C\setminus\{0\}$. Then we have, for $z\in L(N_x)$,
$$\Omega(z+N_x) - \Omega(z-N_x) = \frac{(z+N_x)^j - (z-N_x)^j}{(|z|^2+1)^{j/2}} = N_x^j\,
\frac{\bigl(\frac{z }{N_x} + 1 \big)^j - \big(\frac{z }{N_x}-1\big)^j}{(|z|^2+1)^{j/2}}
$$
Hence, using the change of variable given by the rotation $y = z/(iN_x)$ (which transforms $L_{N_x}$ into $\R$), we obtain
\begin{align*}
C_K(N_x) & = \int_\R N_x^j\,
\frac{\bigl(iy + 1 \big)^j - \big(iy-1\big)^j}{2(y^2+1)^{(j+1)/2}}\,dy
= \frac{(iN_x)^j}2\,\int_\R 
\frac{\bigl(y -i \big)^j - \bigl(y +i \big)^j}{(y-i)^{(j+1)/2} (y+i)^{(j+1)/2}}\,dy\\
& = \frac{(iN_x)^j}2\,\int_\R \Bigg( \frac{\bigl(y -i \big)^{(j-1)/2}}{(y+i)^{(j+1)/2}}
- \frac{\bigl(y +i \big)^{(j-1)/2}}{(y-i)^{(j+1)/2}}\Bigg)
\,dy.
\end{align*}
By the residue theorem, it is easy to check that the last integral equals $2\pi i$, and thus
$$C_K(N_x) = \pi i^{j+1} N_x^j = (-1)^{(j+1)/2} \pi N_x^j.$$
\end{itemize}
\vv

Some of the known proofs of the jump formulas for the Cauchy and Riesz transforms suggest that
one needs to use complex analysis or Clifford analysis to deduce them (see \cite{Tolsa-llibre} or \cite{HMT}, for example), and
$E$ has to coincide with the boundary of a ``reasonable" domain. For example, in \cite{HMT} the authors ask these domains to have locally finite perimeter and their boundaries to be uniformly rectifiable.
The proof   of Theorem \ref{teo1} in the present paper is of real variable nature and avoids such assumptions. I think that
this approach is rather elementary and has its own interest, even in the case when one asks $E$ to be as in \cite{HMT} and the result is not new. 

An important tool for the proof of Theorem \ref{teo1} is the boundedness of the maximal singular integrals
from $M(\R^{n+1})$ to $L^{1,\infty}(\mu)$, which allows an easy reduction of the proof to the case
where $E$ is a subset of a $C^1$ graph. As far as I know, the idea of using the boundedness from 
$M(\R^{n+1})$ to $L^{1,\infty}(\mu)$ in connection with  principal values for singular integrals stems from a paper from
Mattila and Melnikov \cite{MaMe} on the existence of principal values for the Cauchy transform.
See \cite[Chapter 8]{Tolsa-llibre} for a more modern exposition of such techniques.


\section{Proof of Theorem \ref{teo1}}

In the arguments below we allow all the implicit constants in the relation $\lesssim$ to depend on the Calder\'on-Zygmund constants of the kernel and also on the aperture parameter $a$ and on the constant $b$. 

To prove the theorem we can assume that $E$ is a compact subset of a $C^1$ (possibly rotated) graph $\Gamma \subset\R^{n+1}$ with slope at most $1/10$. 
Further, it is enough to prove that for, 
$\HH^n$-a.e.\ $x\in E$,
\begin{equation}\label{eqg1}
\lim_{X_a^+(x)\ni y\to x} \bigl(T_{b|x-y|}\nu(y) + T_{b|x-y|}\nu(2x-y)\bigr) = 2\pv T\nu(x),
\end{equation}
and that
\begin{equation}\label{eqg2}
\lim_{X_a^+(x)\ni y\to x} \bigl(T_{b|x-y|}\nu(y) - T_{b|x-y|}\nu(2x-y)\bigr) = 2C_K(N_x)\,f(x),
\end{equation}
with $C_K(N_x)$ and $f$ as defined in the statement of the theorem. Notice that $2x-y$ is the point symmetric to $y$ with respect to $x$, and thus it belongs to $X_a^-(x)$ if $y\in X_a^+(x)$.
In fact, from the preceding identities one deduces that $T^+\nu(x)\equiv \lim_{X_a^+(x)\ni y\to x} T_{b|x-y|}\nu(y)$
exists for $\HH^n$-a.e. $x\in E$ and that 
\begin{equation}\label{eqg3}
T^+\nu(x) = \pv T\nu(x) + C_K(N_x)\,f(x).
\end{equation}
It is immediate to check also that $T^-\nu(x)\equiv \lim_{X_a^+(x)\ni y\to x} T_{b|x-y|}\nu(2x-y)$ and thus, again
by \rf{eqg1} and \rf{eqg2}, it exists and 
\begin{equation}\label{eqg4}
T^-\nu(x) = \pv T\nu(x) - C_K(N_x)\,f(x).
\end{equation}
Clearly, \rf{eqg3} and \rf{eqg4} are equivalent to \rf{eqga} and \rf{eqgb}.

To prove \rf{eqg1} and \rf{eqg2}, write $\mu=\HH^n|_E$, $y_x^*=2x-y$. For a given $\delta>0$ we consider the maximal operators
$$S_\delta \nu(x) = \sup_{y\in X_a^+(x),|x-y|\leq\delta} \Big| \pv T\nu(x) -\frac12
\bigl(T_{b|x-y|}\nu(y) + T_{b|x-y|}\nu(y_x^*)\bigr)\Big|,$$
$$\wt S_\delta \nu(x) = \sup_{y\in X_a^+(x),|x-y|\leq\delta} \Big| C_K(N_x) \frac{d\nu_{abs}}{d\mu}(x) -\frac12
\bigl(T_{b|x-y|}\nu(y) - T_{b|x-y|}\nu(y_x^*)\bigr)\Big|.$$
By standard arguments, it suffices to show that, for all $\lambda>0$,
\begin{equation}\label{eqlam}
\mu\big(\big\{x\in E: S_\delta\nu(x) + \wt S_\delta\nu(x)>\lambda\big\}\big)\to 0 \quad \mbox{ as $\delta\to0$.}
\end{equation}

Note first that $S_\delta$ and $\wt S_\delta$ are subadditive. Moreover, they are bounded from $M(\R^{n+1})$ to $L^{1,\infty}(\mu)$ uniformly on $\delta>0$. Indeed,  for all $x\in\Gamma$ and $y\in X_a^+(x)$, by standard estimates,
$$|T_{b|x-y|}\nu(y)| + |T_{b|x-y|}\nu(y_x^*)| \lesssim T_*\nu(x) + M_n\nu(x),$$
where $M_n$ is the maximal operator
$M_n\nu(x) = \sup_{r>0}\dfrac{|\nu|(B(x,r))}{r^n}.$
Hence, 
$$S_\delta\nu(x)\lesssim T_*\nu(x) + M_n\nu(x)$$
 for all $x\in\Gamma$. Since both $T_*$ and $M_n$ are bounded from $M(\R^{n+1})$ to $L^{1,\infty}(\mu)$, the same holds for
 $S_\delta$. Regarding $\wt S_\delta$, we argue analogously and take also into account that, by the smoothness of $\Omega(\cdot)$, it holds that $C_K(N_x)$ is uniformly bounded, and thus
$$\Big\| C_K(N_x) \frac{d\nu_{abs}}{d\mu}\Big\|_{L^1(\mu)} \lesssim \|\nu_{abs}\|\leq \|\nu\|.$$

Given an arbitrary $\ve>0$, let $\tau>0$ be such that
$|\nu|(\UU_\tau(E)\setminus E)\leq \ve,$
where $\UU_\tau(E)$ stands for the open $\tau$-neighborhood of $E$. Write
\begin{align*}
\mu\big(\big\{x\in E:& S_\delta\nu(x) + \wt S_\delta\nu(x)>\lambda\big\}\big) \\
&\leq
\mu\big(\big\{x\in E: S_\delta(\chi_{(\UU_\tau(E))^c}\nu)(x) + \wt S_\delta(\chi_{(\UU_\tau(E))^c}\nu) (x)>\lambda/3\big\}\big)\\
& \quad+ \mu\big(\big\{x\in E: S_\delta(\chi_{\UU_\tau(E)\setminus E}\nu)(x) + \wt S_\delta(\chi_{\UU_\tau(E)\setminus E}\nu)(x) >\lambda/3\big\}\big)\\
& \quad+ \mu\big(\big\{x\in E: S_\delta(\chi_{E}\nu)(x) + \wt S_\delta(\chi_{E}\nu)(x)>\lambda/3\big\}\big)\\
& =: A_1(\delta) + A_2(\delta) + A_3(\delta).  
\end{align*}
By the continuity of $T(\chi_{(\UU_\tau(E))^c}\nu)$ in $\UU_{\tau/2}(E)$ and the fact that $|(\chi_{(\UU_\tau(E))^c}\nu_{abs})|(E)=0$, it follows easily that $A_1(\delta)= 0$
for $\delta$ small enough
(depending on $\tau$ and thus on $\ve$).
Also, by the boundedness of $S_\delta$ and $\wt S_\delta$ from $M(\R^{n+1})$ to $L^{1,\infty}(\mu)$,
$$A_2(\delta) \lesssim \frac{|\nu|(\UU_\tau(E)\setminus E)}\lambda \lesssim \frac{\ve}\lambda.$$

To estimate $A_3(\delta)$, let $V$ be an open set such that $\mu(V)<\ve$ and $|\nu_s|(V^c)=0$, and consider a compact
set $K\subset  V\cap E$ such that $|\nu|(V\cap E \setminus K) \leq \ve$. Also, take a
compactly supported $C^1$ function $g$ on $\Gamma$ such that
$$\Bigl\|\chi_E\,\frac{d\nu_{abs}}{d\mu} - g\Big\|_{L^1(\mu)}\leq\ve.$$ 
Then we split
$$\chi_E\nu = g\mu + (\chi_E\nu_{abs} - g\mu) + \chi_{V\cap E\setminus K}\nu_s + \chi_K\nu_s.$$
We have
\begin{align*}
\mu\big(A_3(\delta)\big) &\leq \mu\big(\big\{x\in E: S_\delta(g\mu)(x) + \wt S_\delta(g\mu)(x) >\lambda/12\big\}\big)\\
&\quad + \mu\big(\big\{x\in E: S_\delta(\chi_E\nu_{abs}-g\mu)(x) + \wt S_\delta(\chi_E\nu_{abs}-g\mu) (x)>\lambda/12\big\}\big)\\
& \quad+ \mu\big(\big\{x\in E: S_\delta(\chi_{V\cap E\setminus K}\nu_s)(x) + \wt S_\delta(\chi_{V\cap E\setminus K}\nu_s)(x)>\lambda/12\big\}\big)\\
&\quad + \mu\big(\big\{x\in E\setminus V: S_\delta(\chi_K\nu_s)(x) + \wt S_\delta(\chi_K\nu_s) (x)>\lambda/12\big\}\big) + \mu(V)\\
& =: B_1(\delta) + B_2(\delta) + B_3(\delta) + B_4(\delta)+ \mu(V).
\end{align*}
Again by the boundedness of $S_\delta$ and $\wt S_\delta$ from $M(\R^{n+1})$ to $L^{1,\infty}(\mu)$,
$$B_2(\delta)+B_3(\delta)\lesssim \frac1\lambda \,\Bigl\|\chi_E\,\frac{d\nu_{abs}}{d\mu} - g\Big\|_{L^1(\mu)}+
\frac{|\nu|(V\cap E\setminus K)}{\lambda}\lesssim\frac{\ve}\lambda.$$
Also, by continuity, $B_4(\delta)= 0$ for $\delta$ small enough.

Summarizing, we have shown that
$$\mu\big(\big\{x\in E: S_\delta\nu(x) + \wt S_\delta\nu(x)>\lambda\big\}\big)\lesssim 
\mu\big(\big\{x\in E: S_\delta(g\mu)(x) + \wt S_\delta(g\mu)(x)>\lambda/12\big\}\big) +\ve + \frac\ve\lambda$$
for $\delta$ small enough. We will prove now that, for every $x\in E$ for which $\pv T_\mu g(x)$ exists,
$$S_\delta(g\mu)(x) + \wt S_\delta(g\mu)(x)\to 0\quad \mbox{ as $\delta\to 0$.}$$
This will imply that 
$$\mu\big(\big\{x\in E: S_\delta(g\mu)(x) + \wt S_\delta(g\mu)(x)>\lambda/12\big\}\big)\leq\ve$$
for $\delta$ small enough and will conclude the proof of \rf{eqlam}.
 
First we deal with $S_\delta(g\mu)(x)$.
So fix $x\in E$ such that $\pv T_\mu g(x)$ exists, and given $\ve>0$ let $\delta>0$ be small enough so that
\begin{equation}\label{eqAA}
|\pv T (g\mu)(x) - T_{A\delta}(g\mu)(x)|\leq \ve,
\end{equation}
for some $A\gg1$ to be fixed below. For a fixed $y\in X_a^+(x)$, with $|x-y|\leq\delta$, denote
$ d_y=|x-y|$ and observe that $T_{b d_y}(g\mu)(y)=T(g\mu)(y)$ and $T_{b d_y}(g\mu)(y_x^*) = T(g\mu)(y_x^*)$
for $\delta$ small enough (depending on the continuity of the tangent at $x$), because $g\mu$ is supported on $\Gamma$, which is a $C^1$ graph.
Write
\begin{align*}
\Big| \pv T(g\mu)(x) -\frac12
\bigl(T_{b d_y}(g\mu)(y) + T_{b d_y}\nu(y_x^*)\bigr)\Big| & \leq
\big| \pv T(g\mu)(x) - T_{A d_y}(g\mu)(x) \big| \\&\quad+ 
\frac12\, \big| T(\chi_{\bar B(x,A d_y)^c}g\mu)(x) - T(\chi_{\bar B(x,A d_y)^c}g\mu) (y)\big| \\
&\quad+  
\frac12 \,\big| T(\chi_{\bar B(x,A d_y)^c}g\mu)(x) - T(\chi_{\bar B(x,A d_y)^c}g\mu) (y_x^*)\big|\\
& \quad + \frac12\, \big| T(\chi_{\bar B(x,A d_y)}g\mu)(y) + T(\chi_{\bar B(x,A d_y)}g\mu) (y_x^*)\big|.
\end{align*}
The first term on the right hand side is smaller or equal that $\ve$. The second one 
satisfies
\begin{equation}\label{eqst34}
\big| T(\chi_{\bar B(x,A d_y)^c}g\mu)(x) - T(\chi_{\bar B(x,A d_y)^c}g\mu) (y)\big| 
\lesssim \frac1A\,\|g\|_\infty,
\end{equation}
by standard estimates, taking into account that $|x-y|= d_y\leq \frac1A\dist(x,\supp(\chi_{\bar B(x,A d_y)^c}g\mu))$
and the polynomial growth of $\mu$. Indeed, we just have to write
\begin{align*}
\big| T(\chi_{\bar B(x,A d_y)^c}g\mu)(x) - T(\chi_{\bar B(x,A d_y)^c}g\mu) (y)\big| &
\lesssim \int_{\bar B(x,A d_y)^c} \big|K(x-z) - K(y-z)\big|\,|g(z)|\,d\mu(z)\\&\lesssim \|g\|_\infty \int_{\bar B(x,A d_y)^c}\frac{|x-y|}{|x-z|^{n+1}}\,d\mu(z)
\end{align*}
and estimate the last integral by splitting the domain of integration into annuli centered at $x$, say.
The same estimate holds replacing $y$ by $y_x^*$, and thus choosing $A=\ve^{-1}$ we get
\begin{multline}\label{eqcla89}
\Big| \pv T(g\mu)(x) -\frac12
\bigl(T_{b d_y}(g\mu)(y) + T_{b d_y}\nu(y_x^*)\bigr)\Big| \\
\leq \ve + C\ve\,\|g\|_\infty 
+ \frac12\, \big| T(\chi_{\bar B(x,A d_y)}g\mu)(y) + T(\chi_{\bar B(x,A d_y)}g\mu) (y_x^*)\big|.
\end{multline}

To estimate the last term above we write
\begin{align}\label{eqa46}
\big| T(\chi_{\bar B(x,A d_y)}g\mu)(y) + T(\chi_{\bar B(x,A d_y)}g\mu) (y_x^*)\big| & \leq |g(x)|
\big| T(\chi_{\bar B(x,A d_y)}\mu)(y) + T(\chi_{\bar B(x,A d_y)}\mu) (y_x^*)\big| \\
& \quad+
\big| T(\chi_{\bar B(x,A d_y)}(g-g(x))\mu)(y)\bigr|\nonumber\\
&\quad + \bigl|T(\chi_{\bar B(x,A d_y)}(g-g(x))\mu) (y_x^*)\big|.\nonumber
\end{align}
It is immediate to check that the two last summands on the right hand side are bounded above
by 
\begin{equation}\label{eqa47}
C\!\sup_{z\in \bar B(x,A d_y)}|g(z)-g(x)| \,\frac{\mu(B(x,A d_y))}{ d_y^n}\lesssim A^n\!\sup_{z\in \bar B(x,A\delta)}|g(z)-g(x)|= \ve^{-n}\!\!\sup_{z\in \bar B(x,A\delta)}|g(z)-g(x)|.
\end{equation}
By the continuity of $g$, the right hand side is at most $\ve$ if $\delta$ is small enough. 

Hence it just remains to estimate the first term on the right hand side of \rf{eqa46}. We will compare this term to the analogous one replacing $\mu$ by
$\HH^n|_{L_x}$ (recall that $L_x$ is the tangent $n$-plane of $\Gamma$ at $x$). Notice that the reflection
$R_x:z\mapsto 2x-z$ leaves invariant the measure $\HH^n|_{L_x}$ (i.e., the image measure $R_x\#\HH^n|_{L_x}$ equals 
$\HH^n|_{L_x}$) and also the ball $B(x,A d_y)$. Thus,
\begin{align*}
T(\chi_{\bar B(x,A d_y)}\HH^n|_{L_x}) (y_x^*) & = \int_{B(x,A d_y)} K(y_x^* - z) \,d\HH^n|_{L_x}(z) \\
& =
 \int_{B(x,A d_y)} K(y_x^* - z) \,dR_x\#\HH^n|_{L_x}(z) \\ 
 &= \int_{B(x,A d_y)} K(y_x^* - (2x-z)) \,d\HH^n|_{L_x}(z)\\
 &  = \int_{B(x,A d_y)} K(z-y) \,d\HH^n|_{L_x}(z)
 = -T(\chi_{\bar B(x,A d_y)}\HH^n|_{L_x})(y),
\end{align*}
by the antisymmetry of $K$. Thus,
\begin{align}\label{eqd41}
&\big| T(\chi_{\bar B(x,A d_y)}\mu)(y) + T(\chi_{\bar B(x,A d_y)}\mu) (y_x^*)\big|\\ &\! =
\big| \big(T(\chi_{\bar B(x,A d_y)}\mu)(y) + T(\chi_{\bar B(x,A d_y)}\mu) (y_x^*)\big) - 
\big(T(\chi_{\bar B(x,A d_y)}\HH^n|_{L_x})(y) + T(\chi_{\bar B(x,A d_y)}\HH^n|_{L_x}) (y_x^*)\big)
\big| \nonumber\\
& \!\leq
\big|T(\chi_{\bar B(x,A d_y)}\mu)(y) - 
T(\chi_{\bar B(x,A d_y)}\HH^n|_{L_x})(y) \big| +
\big|T(\chi_{\bar B(x,A d_y)}\mu)(y_x^*) - 
T(\chi_{\bar B(x,A d_y)}\HH^n|_{L_x})(y_x^*) \big|.\nonumber
\end{align}

We will now estimate the first term on the right hand side of \rf{eqd41}. To show that this is small is a routine task.
Since $\Gamma$ is the graph of a $C^1$ function, there is some $\delta>0$ small enough such that, for all $w\in \Gamma\cap B(x,\delta)$, the angle between the tangent planes $L_w$ and $L_x$ is at most $A^{-n-1}\ve$.
As a consequence, $\Gamma\cap B(x,\delta)$ coincides with the graph of a $C^1$ function $F_x:L_x \to \R$ with
slope at most $CA^{-n-1}\ve$. More precisely, denoting by $\Pi$ the orthogonal projection on $L_x$ and writing $\wt F_x(z) = z+ 
F_x(z)N_x$ for $z\in L_x$, we have
$$\Gamma\cap B(x,\delta) = \wt F_x(\Pi(\Gamma\cap B(x,\delta))),$$
with $\|\nabla F_x\|_\infty\lesssim  CA^{-n-1}\ve$.
Then we can write
$$T(\chi_{\bar B(x,A d_y)}\mu)(y) = \int_{\Pi(\bar B(x,A d_y)\cap \Gamma)} K(y-\wt F_x(z))\,d\Pi\#\HH^n|_\Gamma(z).
$$
Therefore, 
\begin{align*}
\big|T(\chi_{\bar B(x,A d_y)}\mu)(y) - 
T(&\chi_{\bar B(x,A d_y)}\HH^n|_{L_x})(y) \big| \\
&  \leq 
\bigg|\int_{\Pi(\bar B(x,A d_y)\cap \Gamma)} \big(K(y-\wt F_x(z)) - K(y-z)\big)\,d\Pi\#(\HH^n|_\Gamma)(z)\bigg|\\
&\quad +\bigg|\int_{\Pi(\bar B(x,A d_y)\cap \Gamma)} K(y-z)\,d(\Pi\#(\HH^n|_\Gamma) -\HH^n|_{L_x})(z) \bigg|\\
& \quad+ \HH^n\big((\bar B(x,A d_y)\cap L_x)\setminus (\Pi(\bar B(x,A d_y)\cap \Gamma))\big)\,\sup_{z\in L_x}|K(y-z)|\\
&=: I_1 + I_2+I_3.
\end{align*}
 Since $|y-z|\gtrsim |y-x|= d_y$ for $z\in L_x$, we derive
$$\sup_{z\in L_x}|K(y-z)|\lesssim \frac1{ d_y^n}.$$
Observe also that, by the fact that the slope of $F_x$ is at most $A^{-n-1}\ve$, for all $w\in \Gamma\cap B(x,A d_y)$,
$$\dist(w, L_x)\lesssim\,A^{-n-1}\ve A d_y  = A^{-n} \ve \, d_y.$$
Then it easily follows that, for some absolute constant $C$
\begin{align*}
\HH^n\big((\bar B(x,A d_y)\cap L_x)\setminus \Pi(&\bar B(x,A d_y)\cap \Gamma)\big) \\
&\leq \HH^n\big(B(x,A d_y)\cap L_x\setminus B(x,A d_y- CA^{-n} \ve \, d_y)\big)
\lesssim \ve\, d_y^n,
\end{align*}
and thus $I_3\lesssim\ve$. Also, taking into account that 
$$d\Pi\#(\HH^n|_\Gamma)(z) = \big(1+|\nabla F_x(z)|^2\big)^{1/2}\,d\HH^n|_{L_x}(z) $$
and again that $\|\nabla F_x\|_\infty\lesssim A^{-n-1}\ve$, we derive
$$I_2\lesssim \sup_{z\in L_x}|K(y-z)|\,A^{-n-1}\ve\,\HH^n(B(x,A d_y)\cap L_x) \lesssim \frac1{ d_y^n}\,\,A^{-n-1}\ve\,(A d_y)^n
 \leq\ve.$$
Concerning $I_1$, we use the fact that for $z\in \Pi(\bar B(x,A d_y)\cap \Gamma)$, 
$$|z-\wt F_x(z)|\lesssim
A^{-n-1}\ve A d_y\leq A^{-n}\ve d_y \ll |y-z|,$$ and thus
$$\big|K(y-\wt F_x(z)) - K(y-z)\big| \lesssim \frac{|\wt F_x(z) - z|}{|y-z|^{n+1}}\lesssim  \frac{A^{-n}\ve d_y}{ d_y^{n+1}}
= \frac{A^{-n}\ve}{ d_y^n}.$$
Hence,
$$I_1\lesssim \frac{A^{-n}\ve}{ d_y^{n}}\,\mu(\bar B(x,A d_y)) \lesssim \ve.$$
So we have shown that
$$\big|T(\chi_{\bar B(x,A d_y)}\mu)(y) - 
T(\chi_{\bar B(x,A d_y)}\HH^n|_{L_x})(y) \big| \lesssim\ve$$
for $\delta$ small enough. The same estimates above work replacing $y$ by $y_x^*$, and thus by \rf{eqd41} it follows that 
$$\big| T(\chi_{\bar B(x,A d_y)}\mu)(y) + T(\chi_{\bar B(x,A d_y)}\mu) (y_x^*)\big|\lesssim\ve.$$
Together with \rf{eqcla89} and \rf{eqa46} (and recalling also that the last two terms in \rf{eqa46} are bounded by $\ve$), we deduce that, for all $y\in X_a(x)\cap B(x,\delta)$,
$$\Big| \pv T(g\mu)(x) -\frac12
\bigl(T_{b d_y}(g\mu)(y) + T_{b d_y}\nu(y_x^*)\bigr)\Big| \leq \ve + C\ve\,\|g\|_\infty$$
for $\delta$ small enough, which proves that $S_\delta(g\mu)(x)\to0$ as $\delta\to0$, as wished.
\vv

Next we have to show that $\wt S_\delta(g\mu)(x)\to 0$ as $\delta\to 0$.
Again we take $A=\ve^{-1}$ and $\delta$ small enough and  we fix
 $y\in X_a^+(x)$, with $|x-y|= d_y\leq\delta$. We want to show that
$$\Big| C_K(N_x) g(x) -\frac12
\bigl(T(g\mu)(y) - T(g\mu)(y_x^*)\bigr)\Big|\leq \ve$$
if $\delta$ is small enough, depending on $\ve$. First we write
\begin{align}\label{eqa91}
\Big| C_K(N_x) g(x)& -\frac12
\bigl(T(g\mu)(y) - T(g\mu)(y_x^*)\bigr)\Big|\\
&\leq \Big| C_K(N_x) g(x) -\frac12
\bigl(T(\chi_{\bar B(x,A d_y)}g\mu)(y) - T(\chi_{\bar B(x,A d_y)}g\mu)(y_x^*)\bigr)\Big| \nonumber\\
& \quad + \frac12\big|T(\chi_{\bar B(x,A d_y)^c}g\mu)(y) - T(\chi_{\bar B(x,A d_y)^c}g\mu)(y_x^*)\big|.
\nonumber
\end{align}
Arguing as in \rf{eqst34} we deduce that
$$\big| T(\chi_{\bar B(x,A d_y)^c}g\mu)(y) - T(\chi_{\bar B(x,A d_y)^c}g\mu) (y_x^*)\big| 
\lesssim \frac1A\,\|g\|_\infty = \ve \,\|g\|_\infty.
$$
Concerning the first term on the right hand side of \rf{eqa91}
we split
\begin{align*}
\Big| C_K(N_x) g(x) -\frac12&
\bigl(T(\chi_{\bar B(x,A d_y)}g\mu)(y) - T(\chi_{\bar B(x,A d_y)}g\mu)(y_x^*)\bigr)\Big|\\
& \leq
\frac12\,\big| T(\chi_{\bar B(x,A d_y)}(g-g(x))\mu)(y)\bigr| + \frac12\,\bigl|T(\chi_{\bar B(x,A d_y)}(g-g(x))\mu) (y_x^*)\big|\\
&\quad+\frac12\,|g(x)|\,\big|T(\chi_{\bar B(x,A d_y)}\mu)(y) - 
T(\chi_{\bar B(x,A d_y)}\HH^n|_{L_x})(y) \big| \\
&\quad + \frac12\,|g(x)|\,
\big|T(\chi_{\bar B(x,A d_y)}\mu)(y_x^*) - 
T(\chi_{\bar B(x,A d_y)}\HH^n|_{L_x})(y_x^*) \big|\\
&\quad + |g(x)| \,
\Big| C_K(N_x) -\frac12
\bigl(T (\chi_{\bar B(x,A d_y)}\HH^n|_{L_x})(y) - T (\chi_{\bar B(x,A d_y)}\HH^n|_{L_x})(y_x^*)\bigr)\Big|
 \\
& =  J_1+ J_2+J_3+J_4+J_5.
\end{align*}
The terms $J_1$ and $J_2$ coincide with the last two terms on the right hand of \rf{eqa46} (modulo
the factor $\frac12$), which are bounded above by $\ve$ if $\delta$ is small enough, as shown in \rf{eqa47}. Also, the terms $J_3$ and $J_4$ coincide with the last two terms in \rf{eqd41} (except for the factor $\frac12|g(x)|$) and, as shown above, they do not exceed $|g(x)|\ve$ if $\delta$ is small enough too.
So altogether we obtain
$$\Big| C_K(N_x) g(x) -\frac12
\bigl(T(g\mu)(y) - T(g\mu)(y_x^*)\bigr)\Big|\leq C\,\ve\,\|g\|_\infty + J_5.$$

It remains to estimate $J_5$. To this end, we write
\begin{align}\label{eqj5} 
J_5 & \leq |g(x)| \,
\Big| C_K(N_x) -\frac12
\bigl(T \HH^n|_{L_x}(y) - T \HH^n|_{L_x}(y_x^*)\bigr)\Big|\\
&\quad + \frac12\,|g(x)| \,
\big|T (\chi_{\bar B(x,A d_y)^c}\HH^n|_{L_x})(y) - T (\chi_{\bar B(x,A d_y)^c}\HH^n|_{L_x})(y_x^*)\big|.
\nonumber
\end{align}
Indeed, above one should understand $T \HH^n|_{L_x}(y) - T \HH^n|_{L_x}(y_x^*)$ and $T (\chi_{\bar B(x,A d_y)^c}\HH^n|_{L_x})(y) - T (\chi_{\bar B(x,A d_y)^c}\HH^n|_{L_x})(y_x^*)$ in  a ``BMO
sense''. As in the estimate \rf{eqst34}, by standard estimates we get
$$\big|T (\chi_{\bar B(x,A d_y)^c}\HH^n|_{L_x})(y) - T (\chi_{\bar B(x,A d_y)^c}\HH^n|_{L_x})(y_x^*)\big|
\lesssim \frac1A = \ve.$$

Finally we will show that the first term on the right hand side of \rf{eqj5} vanishes. First note that if
$\wt y\in\R^{n+1}$ is any point which equals the translation of $y$ along any direction parallel to $L_x$,
then
$$T \HH^n|_{L_x}(y) - T \HH^n|_{L_x}(\wt y) = 0.$$
In other words, the function $T \HH^n|_{L_x}$ (which we should understand in a suitable $BMO$ sense) is invariant by translations parallel to $L_x$. Therefore, when computing 
$T \HH^n|_{L_x}(y) - T \HH^n|_{L_x}(y_x^*)$
we can assume that $y$ and $y_x^*$ are of the
form 
$$y= x + t\,N_x \quad\mbox{ and }\quad y_x^* = x - t\,N_x, \qquad\mbox{for some $t>0$.}$$
Then, using also the fact that $K$ is a convolution kernel
\begin{align*} 
T \HH^n|_{L_x}(y) - T \HH^n|_{L_x}(y_x^*) &= T \HH^n|_{L_x}(x+t\,N_x) - T \HH^n|_{L_x}(x-t\,N_x)\\
& = T \HH^n|_{L(N_x)}(t\,N_x) - T \HH^n|_{L(N_x)}(t\,N_x)
\end{align*}
(recall that $L(N_x)$ is the $n$-plane parallel to $L_x$ through the origin). Taking now into 
account that $K$ is homogeneous of degree $-n$, it follows that
$$T \HH^n|_{L_x}(y) - T \HH^n|_{L_x}(y_x^*) = T \HH^n|_{L(N_x)}(N_x) - T \HH^n|_{L(N_x)}(N_x).$$
By the definition \rf{eqckn}, the right hand equals $C_K(N_x)$ and thus the first term on the right
hand side of \rf{eqj5}  vanishes, as wished.

Gathering the estimates above we deduce that
$$\Big| C_K(N_x) g(x) -\frac12
\bigl(T(\chi_{\bar B(x,A d_y)}g\mu)(y) - T(\chi_{\bar B(x,A d_y)}g\mu)(y_x^*)\bigr)\Big|\lesssim \ve +\ve\,\|g\|_\infty$$
for all $y\in X_a^+(x)\cap B(x,\delta)$ if $\delta$ is small enough, which proves that
$\wt S_\delta(x)\to 0$ as $\delta\to0$ and concludes the proof of the theorem.

\enlargethispage{5mm}


\vvv

\end{document}